\documentclass[12pt]{amsart}%

\usepackage{graphicx}
\usepackage{psfrag}
\usepackage{mathrsfs}

\begin{document}

\newcounter{teorema}
\newtheorem{teor}[teorema]{\sc Theorem}

\newcounter{theorem}[section]
\newtheorem{conj}{\sc Conjecture}
\newcounter{\theconj}[conj]
\newtheorem{defi}[theorem]{\sc Definition}
\newtheorem{lema}[theorem]{\sc Lemma}
\newtheorem{prop}[theorem]{\sc Proposition}
\newtheorem{cor}[theorem]{\sc Corrollary}
\newtheorem{teo}[theorem]{\sc Theorem}
\newtheorem{obs}[theorem]{\sc Remark}
\newtheorem{ques}[theorem]{\sc Question}
\newtheorem{coro}[theorem]{\sc Corollary}
\newtheorem{question}{\sc Question}
\def\bp{\noindent{\it Proof. }}
\def\ep{\noindent{\hfill $\fbox{\,}$}\medskip\newline}
\renewcommand{\theequation}{\arabic{section}.\arabic{equation}}
\renewcommand{\thetheorem}{\arabic{section}.\arabic{theorem}}
\newcommand{\eps}{\varepsilon}
\newcommand{\disp}[1]{\displaystyle{\mathstrut#1}}
\newcommand{\fra}[2]{\displaystyle\frac{\mathstrut#1}{\mathstrut#2}}
\newcommand{\dif}{{\rm Diff}}
\newcommand{\ph}{{\mathcal {PH}}^r_m}
\newcommand{\phn}{{\mathcal {PH}}^r}
\newcommand{\phr}[1]{{\rm Par}^{#1}_m}
\newcommand{\Per}{{\mathcal P}er}
\newcommand{\Z}{\mathbb Z}
\newcommand{\R}{\mathbb R}
\newcommand{\N}{\mathbb N}
\newcommand{\s}{\sigma}
\newcommand{\de}{\em}
\newcommand{\W}{\mathcal W}
\newcommand{\F}{\mathcal F}
\newcommand{\C}{\mathscr C}
\newcommand{\KS}{\mathscr{K\!\!S}}
\def\to{\mathop{\rightarrow}}
\def\ord{\mathop{\rm ord}}
\def\diam{\mathop{\rm diam}}
\def\cc{\mathop{\rm cc}}
\title[SOME RESULTS ON THE INTEGRABILITY OF THE CENTER BUNDLE]
{SOME RESULTS ON THE INTEGRABILITY OF THE CENTER BUNDLE FOR
PARTIALLY HYPERBOLIC DIFFEOMORPHISMS}
\author{F.Rodriguez Hertz}
\author{M.Rodriguez Hertz}
\author{R. Ures}
\thanks{This work was partially supported by FCE 9021, CONICYT-PDT 29/220 and  CONICYT-PDT 54/18 grants}
\address{IMERL-Facultad de Ingenier\'\i a\\ Universidad de la
Rep\'ublica\\ CC 30 Montevideo, Uruguay} \email{frhertz@fing.edu.uy}
\email{jana@fing.edu.uy} \email{ures@fing.edu.uy}
\begin{abstract}
We prove, for $f$ a partially hyperbolic diffeomorphism with center
dimension one, two results about the integrability of its central
bundle. On one side, we show that, if $\Omega(f)=M$ and $\dim(M)=3$,
the absence of periodic points implies its unique integrability. On
the opposite side, we prove that any periodic point $p\in \Per(f)$
of large enough period $N$ has an $f^N$-invariant center manifold
(everywhere tangent to the center bundle).

We also obtain, as a consequence of the last result, that there is
an open and dense subset of $C^1$ robustly transitive and partially
hyperbolic diffeomorphisms with center dimension 1,  such that
either the strong stable or the strong unstable foliation is
minimal. This generalizes a result obtained in \cite{bdu} for
three-dimensional manifolds to any dimension.
\end{abstract}
\maketitle

\section{Introduction}

In this paper we shall consider partially hyperbolic diffeomorphisms
with one dimensional center direction $E^c$. By a partially
hyperbolic diffeomorphism we mean  $f\in\dif(M)$, $M$ a closed
manifold,  admitting a non trivial $Df$-invariant splitting of the
tangent bundle $TM=E^s\oplus E^c\oplus E^u$, such that all unit
vectors $v^\s\in E^\s_x$ with $\s=s,c,u$ and $x\in M$ verify:
$$\|Df(x)v^s\|<\|Df(x)v^c\|<\|Df(x)v^u\|$$
for some suitable Riemannian metric, which we call {\de adapted}. It
is also required that the norm of the operators $Df(x)|_{E^s}$ and
$Df^{-1}(x)|_{E^u}$ be strictly less than $1$.  We shall denote
$\phn (M)$ the family of $C^r$  partially hyperbolic diffeomorphisms
of $M$. Along this paper we will consider only the case $\dim E^c=1$
and we denote the set of such diffeomorphisms by $\phn_1 (M)$.

On one hand, it is well known by classical invariant manifold theory
that the bundles $E^s$ and $E^u$ are uniquely integrable thus
obtaining two foliations called the strong stable and the strong
unstable foliations. On the other hand,  it is not known in general
whether either the center bundle, the center stable
($E^{cs}=E^s\oplus E^c$) or the center unstable ($E^{cu}=E^u\oplus
E^c$) are integrable. The hypothesis of integrability of this
bundles has played an important role in partial hyperbolicity
theory, see for instance \cite{hhu2}. Although recent work shows
that the integrability assumption can be bypassed to obtain
ergodicity (see \cite{bw2}, \cite{hhu}) it seems that it remains to
play a crucial role if one looks for a topological description or
even classification of partially hyperbolic diffeomorphisms (see for
instance \cite{brin.burago.ivanov}).

In this paper we prove two results about the integrability of the
center bundle. In our first theorem we prove that, if $\dim(M)=3$
and the nonwandering set is the whole manifold, the absence of
periodic points implies its unique integrability and, in the second,
that periodic points of period $N$ high enough have central curves
(tangent at {\em every} point to $E^c$) invariant by $f^N$.

\begin{teor}\label{nopp} Let $f\in \phn_1(M)$ be such that $\Per (f)=\emptyset$ and $\Omega (f)=M$ and assume
that $\dim (M)=3$. Then, $E^c$ is uniquely integrable.
\end{teor}


\begin{teor} \label{pp} Let $f\in \phn_1$ .
 There exists $K>0$ such
that for any $p \in \Per (f)$ with period $N>K$ there exists,
through $p$, an $f^N$ invariant curve tangent to $E^c$ at every
point.
\end{teor}

A $C^r,\, r\geq 1,$ robustly transitive diffeomorphism is a
diffeomorphism having a neighborhood in $\dif^r(M)$ such that every
$g$ in this neighborhood is transitive. We show as a consequence of
Theorem \ref{pp} a generalization of a result in \cite{bdu} to $M$
of any dimension. In the cited paper the same result is proved for
any $f\in \phn_1(M)$ and $\dim(M)=3$ or for $M$ of any dimension but
assuming unique integrability of the center bundle.

\begin{teor}\label{minimal} Let $\mathscr T (M)$ be the set of $C^1 $ robustly transitive diffeomorphisms.
Then, there exists an open and dense subset of $\mathscr T (M)\cap
\phn_1(M)$ such that either the strong stable or the strong unstable
foliation is minimal.

\end{teor}

\emph{Acknowledgements.} The authors want to thank the Fields
Institute for warm hospitality and financial support during their
visit in January 2006.

\section{Preliminaries}

\par It is a known fact that, for $f\in \phn(M)$, there are
foliations ${\W}^\s$ tangent to the distributions $E^\s$ for
$\s=s,u$ (see for instance \cite{bp}).

Due to Peano's Theorem, for each $x\in M$ there are curves
$\alpha_x(t)$ such that $\alpha_x(0)=x$ and $\dot\alpha_x(t)\in
E^c(\alpha_x(t))\setminus\{0\}$ for some open interval of parameters
$t$ containing $0$. We shall call these curves {\de central curves
through $x$}, and denote by $W^c_{loc}(x)$ the component of a
central curve through $x$ intersected by a small ball. It is easy to
see that $f$ takes central curves into
central curves.\par %
Denoting the leaf of $\W^\s$ through $x$ by $W^\s(x)$, with
$\s=s,u$, we write, as usual, $W^\s_{loc}(x)$ for the connected
component of $W^\s(x)\cap B(x)$, where $B(x)$ is a small ball around
$x$. Observe that for any choice of $W^c_{loc}(x)$, the sets
$$W^\s_{loc}(W^c_{loc}(x))=\bigcup_{y\in W^c_{loc}(x)}W^\s_{loc}(y)\qquad \s=s,u$$
are $C^1$ (local) manifolds tangent to the bundle
$E^{c\s}=E^\s\oplus E^c$ (with $\s=s,u$) at every point (see, for
instance \cite{brin.burago.ivanov}). For further use we will call,
respectively, $W^{cs}_{loc}(x)$ and $W^{cu}_{loc}(x)$ the sets
obtained as above depending, as it is obvious, on the choice of
$W^c_{loc}(x)$.
\begin{obs}\label{remark central en cs}
Moreover, given $x,y\in M$, for all $W^{cs}_{loc}(x)$ such that
$y\in W^{cs}_{loc}(x)$, there exists a central curve $W^c_{loc}(y)$
through $y$ contained in $W^{cs}_{loc}(x)$ (see
\cite{brin.burago.ivanov})
\end{obs}

Let us say that a set $\Gamma$ is {\de $\s$-saturated} if $\Gamma$
is union of leaves of ${\mathcal W}^\s$, $\s=s,u$ and let as call
the accessibility class of $x$, $AC(x)$,  the minimal $s$- and $u$-
saturated set that contains the point $x$ (that is, the  set of
points that can be joined to $x$ by a $us$-path). If $f$ has only
one accessibility class we say that it satisfies the {\de
accessibility property}.

\section{Absence of periodic points}

\begin{lema} \label{nopp.acces.int} Let $f\in \phn_1(M)$ be such that $\Per (f)=\emptyset$ and $\Omega (f)=M$.
Then, either $f$ has the accessibility property or $E^s$ and $E^u$
are jointly integrable.
\end{lema}

\bp Let $\Gamma(f)$ be the set of points such that its accessibility
class is not open and suppose that $\emptyset\neq \Gamma(f)\neq M$.
Thus, Lemma A.5.1 of \cite{hhu} implies the existence of a periodic
point in $\Gamma(f)$ contradicting that $\Per (f)=\emptyset$.

$\Gamma(f)=M$ is equivalent to the joint integrability of $E^s$ and
$E^u$ (see \cite{hhu}) \ep

\begin{obs} \label{clas.densas} Observe that the same proof gives
that, for $f\in \phn_1(M)$ such that $\Per (f)=\emptyset$ and
$\Omega (f)=M$, every closed invariant $su$-saturated set is either
empty or the whole $M$.
\end{obs}

The following lemma generalizes (with essentially the same proof)
Brin's result (\cite{b}) stating that accessibility implies
transitivity. For the sake of completeness we include the proof
here.

\begin{lema} \label{trans} Let $f\in\phn(M)$ be such that $\Omega (f)=M$ and assume that
every closed invariant $su$-saturated set is either empty or the
whole $M$. Then $f$ is transitive.
\end{lema}

\bp Let $U$ and $V$ be two open sets. For all $x\in M$ the set
$K=\cap^\infty_{i=0}\overline{\cup_{n\geq i}AC(f^n(x))}$ is
invariant, closed, $su$-saturated and nonempty (observe that
$AC(f^n(x))=f^n(AC(x))$) and so, $K=M$. Then, by taking $x\in U$, we
can chose $N\in \mathbb{N}$ such that $AC(f^N(x))\cap V\neq
\emptyset$. Call $U_N=f^N(U)$.

We shall show that there exists $n\in \mathbb{N}$ such that
$f^n(U_N)\cap V\neq \emptyset$ which implies $f^{N+n}(U)\cap V\neq
\emptyset$.  Since $U$ and $V$ are arbitrary open sets the
transitivity of $f$ follows from this last property.

The considerations above imply that there is an $su$-path
$[z_0,\dots,z_k]$ with $z_0 \in U_N$ and $z_k \in V$. By continuity
of the strong stable and unstable foliations we can choose $V_0,\,
V_1,\,\dots,\, V_k$ open sets such that:
\begin{itemize}
\item $z_i\in V_i$ $\forall i=0,\,\dots,\,k$.
\item $V_0\subset U_N$ and $V_k \subset V$
\item for each point of $x\in V_i$ there exists a $su$-path $[x=x_i, x_{i+1},\dots,x_k]$ joining
$ x$ with a point of $V_k$ with $x_j\in V_j$ $\forall
j=i,\,\dots,\,k$.
\end{itemize}

Suppose that the path $[z_0,z_1]$ is tangent to the stable bundle
(the unstable case is a little bit easier), then there exists a
neighborhood $B\subset V_1$ of $z_1$ such that each point in it can
be joined with $V_0$ by an $s$-path of uniformly bounded length (in
fact its length can be chosen approximately of the length of
$[z_0,z_1]$) and we can chose $B$ in such a way that there exists
$\rho>0$ such that $W^s(x)\cap V_1\supset W^s_{\rho}(x)$ $\forall
x\in B$. Since $\Omega(f)=M$, there exists an arbitrarily large
$m\in \mathbb{N}$ and a point $w\in B$ such that $f^{-m}(w)\in B$.
This implies that, if $m$ is large enough, $f^{-m}(W^s(w))$ contains
the path joining $f^{-m}(w)$ and $V_0$. Thus $f^{-m}(V_1)\cap
U_N\neq \emptyset$ which implies $f^m(U_N)\cap V_1\neq \emptyset$.

Now substitute $V_1$ by $V_1\cap f^m(U_N)$ and repeat the procedure.
By induction we obtain that there is  $n\in \mathbb{N}$ such that
$f^n(U_N)\cap V\neq \emptyset$.

  \ep

\begin{obs} It is not the issue of this work to achieve the minimal
hypothesis to obtain transitivity by using Brin's argument. However,
let us mention that almost the same proof works if one substitutes
the hypothesis on the density of every  invariant saturated nonempty
set by the weaker one demanding the existence of an accessibility
class whose orbit by $f$ is dense ( there exists $x$ such that
$\overline{\cup\{f^n(AC(x));n\in \mathbb{Z}\}}=M$).
\end{obs}

The following theorem is a direct corollary of Remark
\ref{clas.densas} and Lemma \ref{trans}.

\begin{teo} \label{per.trans} Let $f\in \phn_1(M)$ be such that $\Per (f)=\emptyset$ and $\Omega
(f)=M$. Then, $f$ is transitive.
\end{teo}

Before proving Theorem \ref{nopp} let us state the following lemma,
which is a consequence of continuity and transversality of the
invariant bundles:
\begin{lema}\label{lema integrable}
For $\eps>0$ there exists $\delta>0$ such that if $d(x,y)<\delta$
and $z\in W^c_\delta (x)$, then $W^c_{loc}(y)\cap
W^s_\eps(W^u_\eps(z))\ne\emptyset$, regardless of the choice of
center leaves for $x$ and $y$. \par In particular, if
$W^c_{loc}(y)\subset W^{cu}_{loc}(x)$ then $W^c_{loc}(y)\cap
\W^u_\eps(z)\ne\emptyset$ for all $z\in W^c_\delta(x)$
\end{lema}
\begin{obs} \label{centralcrece} As a corollary of lemma above, if $E^c$, restricted to some $W^{cu}_{loc}(x)$, is non
uniquely integrable at $x$, then for sufficiently small $\delta>0$,
and for each connected central subsegment containing $x$, say $c$,
in one of the two separatrix, there is $N>0$ for which
$f^n(c)\not\subset B_\delta(f^n(x))$ for all $n\geq N$.\end{obs}

\noindent{\it Proof of Remark \ref{centralcrece}.}\quad Take $c_1$
and $c_2$ two different center curves, contained in the same
component of $W^{cu}_{loc}(x)\setminus W^{u}_{loc}(x)$  and having
$x$ as endpoint. Since $c_1$ and $c_2$ are different, there exist
$y_1\neq y_2$ such that $y_2\in W^u_{loc}(y_1)$ and  $y_i\in c_i,\,
i=1,2$. The exponential growth of $W^u(y_1)$ under the action of $f$
implies that there exists $N>0$ such that $f^n(y_2)\notin W^u_\eps
(f^n(y_1))$ for all $n\geq N$ and, by lemma above, we obtain that
$f^n(c_2)\not\subset B_\delta(f^n(x))$.\ep
\medskip
Observe that one can prove without using Lemma \ref{lema integrable}
that either $c_1$ or $c_2$ should grow but, in fact, what is proved
in Remark \ref{centralcrece} is that both center curves grow.

 \noindent{\it Proof of Theorem \ref{nopp}.}\quad
 By Theorem
\ref{per.trans} we know that  $f$ is transitive.

As unique integrability is a local property we can suppose, by
taking a double covering and $f^2$ if necessary, that $E^c$ is
oriented and its orientation is preserved by $f$.

Suppose that $E^c$ is not uniquely integrable at $x$. Then, there
are two different arcs $\alpha$ and $\beta $ tangent to $E^c$
beginning at $x$ with the same (positive) orientation. By taking
intersections of $W^{cs}_{loc}(\alpha)$ and $W^{cs}_{loc}(\beta)$
with some $W^{cu}_{loc}(x)$ we can assume that both arcs are
contained in the same $W^{cu}_{loc}(x)$ (in case
$W^{cs}_{loc}(\alpha)\subset W^{cs}_{loc}(\beta)$ we can do the same
argument for $f^{-1})$

Since $W^{cu}_{loc}(x)$ is two dimensional $\alpha$, $\beta$ and a
conveniently chosen unstable arc bound an open region $U$ of
$W^{cu}_{loc}(x)$ and, for every point $z$ in $U$, there is a center
arc $\gamma$ inside $W^{cu}_{loc}(x)$ joining (in the positive
orientation) $x$ with $z$. As $f$ is transitive, taking the
intersection of the strong stable manifold of a point with dense
forward orbit with $U$, we can choose $z$ such that its forward
orbit is dense and $\gamma$ with length much less than the $\delta$
of Remark \ref{centralcrece}.

Now Remark \ref{centralcrece} implies that for $n>N$ the length of
$f^n(\gamma)$ is larger than $\delta$.  We can take from $z$  a
center continuation of $\gamma$ in the positive direction and a
point $w$ in it, close to $z$,  but not in $\gamma$.

Consider $C=W^u_\varepsilon(\gamma)$ for $\varepsilon$ small.

Then there is $K$ very large (in particular larger that $N$) such
that $f^K(z)$ is very close to $w$. Since the length of
$f^K(\gamma)$ is larger than $\delta$ and the unstable manifolds
growth exponentially, the projection of $f^K(C)$ to
$W^{cu}_{loc}(x)$ contains $C$. This implies that there is a
$f^K$-invariant strong stable manifold and thus, we obtain a
periodic point.

\ep
\begin{obs} In fact with the same argument can be proved that if $f\in
\phn_1(M)$ satisfies $\dim (E^s)=1$, $\Omega (f)=M$ and $\Per
(f)=\emptyset$, $E^{cu}$ is uniquely integrable.
\end{obs}

If $E^s$ and $E^u$ are jointly integrable the  assumptions  on the
dimension of the the strong bundles and the nonwandering set are not
needed to obtain the unique integrability of the one dimensional
center bundle $E^{c}$.

 In order to prove next theorem we need the following standard
lemma:
\begin{lema}\label{lema creacion periodicos}
There is $\eps_0>0$ such that if $x\in\Gamma(f)$ verifies
$f^k(B^{su}_{\eps_0}(x))\cap B^{su}_{\eps_0}(x)\neq\emptyset$ for
some $k>0$, then there is a periodic point in $B^{su}_{\eps_0}(x)$.
\end{lema}

\begin{teo} \label{jointly.integrable}
Let $f\in \phn_1(M)$ be such that $\Per (f)=\emptyset$ and  $E^s$
and $E^u$ are jointly integrable. Then $E^c$ is uniquely integrable.
\end{teo}
\bp If $E^c$ is not uniquely integrable at $x\in M$ then there exist
two central curves $\alpha$ and $\beta$ through $x$. As in the proof
of Theorem \ref{nopp}, by possibly taking intersections, we may
assume, for instance that
$\beta\subset W^u_{loc}(\alpha)$.\par%

Consider three different points $w_2<w_1<w_3$ in $\omega(x)$, such
that $d(w_i, w_j)<\delta/4$, $i,j=1,2,3$.  This is possible since
$f$ has no periodic points, and $E^s\oplus E^u$ is a codimension one
bundle, so we can suppose that $w_1$ is locally between
$W^{su}(w_2)$ and $W^{su}(w_3)$, integral manifolds of $E^s\oplus
E^u$. Take $n_1$ such that $d(f^{n_1}(x),w_1)<\delta/16$ ($\delta$
as in Remark \ref{centralcrece}) We are assuming that $E^c$ is
oriented and that $f$ preserves its orientation (modulo taking a
double covering and $f^2$, if necessary). Take a small arc $\gamma$
in $f^{n_1}(\alpha)$ beginning at $f^{n_1}(x)$. As $E^c$ is not
uniquely integrable at $f^{n_1}(x)$ the previous observation implies
that, we may take $n_2$ large enough so that $f^{n_2}(x)$ be
$\delta/16$-near $w_2$, and the length of $f^{n_2-n_1}(\gamma)$ be
greater than $\delta$. Now, projecting locally via the
$su$-foliation we obtain a map of the interval and as a consequence
there is a point $e$ in an $su$-leaf such that $f^{n_2-n_1}(e)$ is
in the same $su$-leaf. This can be made so that Lemma \ref{lema
creacion periodicos} applies, whence we would
have a periodic point. \par%
This implies the unique integrability at $f^{n_1}(x)$ and, of
course, at $x$ on one direction of $E^c$. To obtain the unique
integrability on the other direction we argue in the same way with
$w_3$ instead  of $w_2$. \ep

\medskip

We remark that in Theorems \ref{nopp} and \ref{jointly.integrable}
we prove not only the existence of a foliation tangent to $E^c$ but
its uniqueness.

\section{Existence of central curves for periodic points}

This section is devoted  to prove Theorem \ref{pp}.
\medskip

\noindent {\em Proof of Theorem \ref{pp}.} \quad By the standard
Center Manifold Theorem, through any periodic point $p$ there exists
an immersed, invariant by the period $N$ of $p$, curve  $\gamma$
such that it is tangent at $p$ to $E^c(p)$ and it is invariant when
it make sense.
 After that, take a connected component
of $\gamma\setminus p$ say $\gamma_1$. Suppose for a while that the
center eigenvalue at $p$ is positive. Then, we have to situations:
either $f^N(\gamma_1)\subset \gamma_1$ or $f^{-N}(\gamma_1)\subset
\gamma_1$ (the simultaneous occurrence of both situations is
possible). Suppose that we are in the first case. Then, as in
\cite{brin.burago.ivanov}, $W^s_{loc}(\gamma_1)=\cup_{x\in \gamma_1}
W^s_{loc}(x)$ is a $C^1$ $f^N$-invariant manifold tangent to
$E^{cs}$ at every point. Analogously, for the second case and
$W^u_{loc}(\gamma_1)$. As a conclusion we have that associated to
each component of $\gamma$ we obtain either a center stable or a
center unstable manifold invariant by $f^N$ or $f^{-N}$
respectively. In case the center eigenvalue were negative we can do
the same procedure  for $f^{2N}$ obtaining that either
$W^s_{loc}(\gamma)$ or $W^u_{loc}(\gamma)$ is respectively $f^N$ or
$f^{-N}$-invariant. In fact, if for example  $W^s_{loc}(\gamma_1)$
is invariant by $f^{2N}$ then  $W^s_{loc}(\gamma_1)\cup
f^N(W^s_{loc}(\gamma_1))$ is invariant by $f^N$.

Suppose that, without loss of generality, $W^s_{loc}(\gamma_1)$ is
$f^N$-invariant.

We shall make use of the following property: given $\varepsilon>0$
there exist $\rho>0$ and $\tau>0$ such that if $y\in W^c_{\rho}(x)$
for some center manifold of $x$, any center curve beginning at
$W^s_{\tau}(y)$ and contained in $W^{cs}(x)$ intersects (with length
near $\rho$) $W^s_{\varepsilon /2}(x)$. Moreover, given $\alpha>0$
there exists $\varepsilon>0$ such that if $x\in W^{cs}_{\alpha}(y)$
then $W^s_{\varepsilon}(x)$ intersects any center manifold $W^c(y)$
contained in $W^{cs}(y)$.

Fix a small $\alpha>0$ and take $\varepsilon$ with the preceding
property and $\rho<\alpha$ and $\tau$ as above. Finally, take $K$
such that $\lambda^K\varepsilon<\tau /3$, where $\lambda$ is such
that the norm of $Df(x)|_{E^s}$ is less than $\lambda<1$.

Now suppose that $N>K$ and take $x\in W^c_r(p)$, $x\neq p$, and
$r<\alpha$ so small that the length of $f^N(W^c_r(p))$ is smaller
than $\rho$.

Then,
\begin{itemize}
\item $W^s_{\varepsilon}(x)$ intersects any center curve beginning at
$p$.

\item $f^N(W^s_{\varepsilon}(x))$ cuts a center manifold that
contains $x$ in a point $y$

\item $f^N(W^s_{\varepsilon}(x))\subset W^s_{\tau}(y)$
\end{itemize}

Let $\{X_n\}$ be a sequence of $C^1$ line fields defined in
$W^{cs}_{loc}(p)$ and converging in the $C^0$ topology to $E^c$.
Observe that, for $n$ big enough, the integral curves of $X_n$
intersecting $W^s_{\tau}(y)$ also cut $W^s_{\varepsilon}(x)$. Then,
consider the following maps
$\varphi_n:W^s_{\varepsilon}(x)\rightarrow W^s_{\varepsilon}(x)$.
First, for $z\in W^s_{\varepsilon}(x)$, take $f^N(z)$ and, after
this, since $f^N(z)\in W^s_{\tau}(y)$, take the point of
intersection of the solution of $X_n$ through $f^N(z)$ with
$W^s_{\varepsilon}(x)$. By Brower's Theorem $\varphi_n$ has a fixed
point. This means that there exists $w_n \in W^s_{\varepsilon}(x)$
such that $f^N(w_n)$ and $w_n$ are in the same integral curve for
$X_n$. Arzela-Ascoli's Lemma implies that we have a limit center
curve and a point $w$ in it such that $f^N(w)\in W^c_{loc}(w)$. If
$f^{kN}(w)\rightarrow_{k\rightarrow \infty} p$ we obtain the
invariant center curve through $p$\, by iteration of the obtained
above. If $f^{kN}(w)$ does not converge to $p$ then, it is easy to
prove, that there is another periodic point $p_1$ in $\gamma_1$ and
a center arc joining  $p$ and $p_1$ that verifies the theorem.\ep

\section{Minimality of strong foliations}

 The proof of Theorem \ref{minimal} is identical to that of the corresponding Theorem of
 \cite{bdu} by
 observing that Theorem \ref{pp} substitutes Lemma 5.2 of \cite{bdu}.

All the discussion in \cite{bdu} about the non orientability of the
center bundle is nowadays solved thanks to the open and denseness of
the diffeomorphisms with the accessibility property among the ones
with one dimensional center bundle (see \cite{dw}, \cite{hhu})


\end{document}